\documentclass[12pt,a4paper]{article}
\usepackage{amsmath,amsxtra,amsfonts,amssymb,latexsym,amsthm,amscd,verbatim,amstext}
\usepackage[mathscr]{eucal}

\newcommand{\dif}{\mathrm{d}}

\input cyracc.def

\newfont{\tricyr}{wncyr10 at 12pt}
\newfont{\tricyi}{wncyi10 at 12pt}
\newfont{\tricyb}{wncyb10 at 12pt}
\newfont{\Tricyr}{wncyr10 at 13.6pt}
\newfont{\Tricyi}{wncyi10 at 13.6pt}
\newfont{\Tricyb}{wncyb10 at 13.6pt}
\newfont{\tricmr}{cmr10 at 13.6pt}
\newfont{\tricmi}{cmti10 at 13.6pt}
\newfont{\tricmb}{cmb10 at 13.6pt}

\theoremstyle{plain}
\newtheorem{Th}{Theorem}

\newtheorem{Lem}{Lemma}

\theoremstyle{definition}

\newtheorem{Not}{Remark}

\begin{document}
\centerline {{\bf The existence and space-time decay rates of strong solutions to }} 
\centerline {{\bf Navier-Stokes Equations in weighed $L^\infty(|x|^\gamma{\rm dx})\cap L^\infty(|x|^\beta{\rm dx})$ spaces}} 

\vskip 0.7cm
\begin{center}
D. Q. Khai, N. M. Tri
\end{center}
\begin{center}
Institute of Mathematics, VAST\\
18 Hoang Quoc Viet, 10307 Cau Giay, Hanoi, Vietnam
\end{center}
\vskip 0.7cm 
{\bf Abstract}: In this paper, we prove some results on the 
existence and \linebreak space-time decay rates of global strong 
solutions of the Cauchy problem for the
Navier-Stokes equations in weighed 
$L^\infty(\mathbb R^d,|x|^\gamma{\rm dx})\cap 
L^\infty(\mathbb R^d,|x|^\beta{\rm dx})$ spaces.
\footnotetext[3]{{\it e-mail address}: Khaitoantin@gmail.com, Triminh@math.ac.vn}
\vskip 0,5cm 
\centerline{\bf \S1.  Introduction} 
\vskip 0,5cm 
This paper studies the Cauchy problem of the incompressible Navier–Stokes equations (NSE)
in the whole space $\mathbb R^d$ for $d\geq 2$,
\begin{align} 
\left\{\begin{array}{ll} \partial _tu  = 
\Delta u - \nabla .(u \otimes u) - \nabla p , & \\ 
\nabla .u = 0, & \\
u(0, x) = u_0, 
\end{array}\right . \notag
\end{align}
which is a condensed writing for
\begin{align} 
\left\{\begin{array}{ll} 1 \leq k \leq d, \ \  \partial _tu_k  
= \Delta u_k - \sum_{l =1}^{d}\partial_l(u_lu_k) - \partial_kp , & \\ 
\sum_{l =1}^{d}\partial_lu_l = 0, & \\
1 \leq k \leq d, \ \ u_k(0, x) = u_{0k} .
\end{array}\right . \notag
\end{align}
The unknown quantities are the velocity 
$u(t, x)=(u_1(t, x),\dots,u_d(t, x))$ of the fluid 
element at time $t$ and position $x$ and the pressure $p(t, x)$.\\
There is an extensive literature on the existence and decay rate of strong  solutions of the Cauchy problem for NSE. 
Maria E. Schonbek \cite{M. Schonbek} established the decay of the homogeneous
 $H^m$ norms for solutions to NSE in two dimensions. She showed that if $u$ is a solution to NSE with an arbitrary  $u_0 \in H^m \cap L^1(\mathbb R^2)$ with $ m \geq 3$ then 
$$
\|D^\alpha u\|^2_2 \leq C_\alpha (t+1)^{-(|\alpha|+1)}\ {\rm and}\ \|D^\alpha u\|_\infty 
\leq C_\alpha (t+1)^{-(|\alpha|+\frac{1}{2})}\ {\rm for\ all}\  t\geq 1, \alpha \leq m.
$$
Zhi-Min Chen \cite{Z. Chen} showed that if $u_0 \in L^1(\mathbb R^d)\cap L^p(\mathbb R^d), (d\leq p<\infty)$ 
and  $\|u_0\|_1 + \|u_0\|_p$ is small enough then there is a unique solution \linebreak
 $u \in BC([0, \infty);L^1 \cap L^p)$, which satisfies decay property
$$
\underset{t >0}{\rm sup}t^{\frac{d}{2}}\big(\|u\|_\infty +t^{\frac{1}{2}}\|Du\|_\infty + t^{\frac{1}{2}}\|D^2u\|_\infty \big) < \infty.
$$
Kato \cite{T. Kato 1984}  studied strong solutions in the spaces $L^q(\mathbb R^d)$ by applying the $L^q - L^p$ 
estimates for the semigroup generated by the Stokes operator. He showed that there is $T>0$ and a unique
solution $u$, which satisfies 
\begin{gather*}
t^{\frac{1}{2}(1-\frac{d}{q})}u \in BC([0,T); L^q), \ {\rm for}\  d \leq q \leq \infty, \\
t^{\frac{1}{2}(2-\frac{d}{q})}\nabla u \in BC([0,T); L^q), \ {\rm for}\  d \leq q \leq  \infty,
\end{gather*}
as $u_0 \in L^d(\mathbb R^d)$. He showed that $T = \infty$ if $\big\|u_0\big\|_{L^d(\mathbb R^d)}$ is small enough.\\
In 2002, Cheng He and Ling Hsiao \cite{Cheng He} extended the results of Kato, 
they estimated on decay rates of higher order derivatives in time variable and space
 variables for the strong solution to NSE with initial data in $L^d(\mathbb R^d)$. They showed that if $\big\|u_0\big\|_{L^d(\mathbb R^d)}$ 
is small enough then there is a unique solution $u$, which satisfies 
\begin{gather*}
t^{\frac{1}{2}(1+ |\alpha| + 2\alpha_0-\frac{d}{q})} D_x^\alpha D^{\alpha_0}_tu \in BC([0, \infty); L^q), \ {\rm for}\  q \geq d, \\
t^{\frac{1}{2}(2+|\alpha|-\frac{d}{q})}D_x^\alpha p \in BC([0, \infty); L^q), \ {\rm for}\   q \geq d,
\end{gather*}
where $\alpha = (\alpha_1, \alpha_2,...,\alpha_d), |\alpha| = \alpha_1 + \alpha_2+...+ \alpha_d$ and $\alpha_0 \in \mathbb N$. $D^\alpha_x$ 
denotes $\partial_x^{|\alpha|}= \partial^{|\alpha|}/(\partial_{x_1}^{\alpha_1}\partial_{x_2}^{\alpha_2}...\partial_{x_d}^{\alpha_d}),\ 
\partial_t^{\alpha_0}= \partial^{\alpha_0}/\partial t^{\alpha_0}$.\\
In 2005,  Okihiro Sawada \cite{O. Sawada} obtained the decay rate of solution to NSE
 with initial data in $\dot{H}^{\frac{d}{2}-1}(\mathbb R^d)$. 
He showed that every mild solution in the class
$$
u \in BC([0, T); \dot{H}^{\frac{d}{2}-1})\ {\rm and}\ t^{\frac{1}{2}(\frac{d}{2}-\frac{d}{p})} u 
\in BC([0, T); \dot{H}^{\frac{d}{2}-1}_p), 
$$
for some $T>0$ and $p\in (2,\infty]$ satisfies
$$
\big\|u(t)\big\|_{\dot{H}^{\frac{d}{2}-1+\alpha}_q} \leq K_1(K_2\tilde \alpha)^{\tilde\alpha}t^{-\frac{\tilde \alpha}{2}}\ {\rm for}\ \alpha >0, q\geq 2,\ {\rm and}\ \tilde \alpha: =\alpha +\frac{d}{2}-\frac{d}{q} 
$$
where constants $K_1$ and $K_2$ depend only on $d,p, M_1$, and $M_2$ with \linebreak $M_1 =  \underset{0<t<T}{\rm sup}\big\|u(t)\big\|_{\dot{H}^{\frac{d}{2}-1}}$ and $M_2 =  \underset{0<t<T}{\rm sup}t^{\frac{d}{2}(\frac{1}{2}-\frac{1}{p})}\big\|u(t)\big\|_{\dot{H}^{\frac{d}{2}-1}_p}$.\\
The time-decay properties are therefore well understood. However, there are few results
on the spatial decay properties. Farwing and Sohr \cite{R. Farwig} showed a class of weighted
$|x|^\alpha$ weak solutions with second derivatives in space variables and one order
derivatives in time variable in $L^s([0,+\infty);L^q)$ for $1 < q < 3/2, 1 < s < 2$ and
$0\leq 3/q +2/s-4\leq \alpha <{\rm min}\{1/2,3-3/q\}$ in the case of exterior domains. 
In \cite{R. Farwig and H. Sohr},
they also showed that there exists a class of weak solutions satisfying 
\begin{align} 
\big\||x|^{\frac{\alpha}{2}}u\big\|^2_2 +\int^t_0\big\||x|^{\frac{\alpha}{2}}\nabla u\big\|^2_2{\rm d}t \leq
\left\{\begin{array}{ll} C(u_0,f,\alpha)  \ \ \  \ \ \ \ \  \ \ \  \ \ \ \ \ \ \ \ \ 
\ {\rm if} \ 0\leq \alpha <\frac{1}{2}, & \\ 
C(u_0,f,\alpha',\alpha)t^{\frac{\alpha'}{2}-1/4}\ \  {\rm if}
\  \frac{1}{2} \leq \alpha <\alpha' < 1,  & \\
C(u_0,f)(t^{1/4}+t^{1/2}) \ \ \ \ \ \ \ \ \ \ \ \ \ \  \ {\rm if}\ \alpha =1.
\end{array}\right . \notag
\end{align}
While in \cite{C. He}, a class of weak solutions 
$$
(1+|x|^2)^{1/4}u \in L^{\infty}([0,+\infty); L^p(\mathbb R^3))
$$
was constructed for $6/5 \leq p < 3/2$, which satisfies 
$$
\big\||x|^{\frac{1}{2}}u\big\|^2_2 +\int^t_0\big\||x|^{\frac{1}{2}}
\nabla u\big\|^2_2{\rm d}t \leq C(u_0,f)(t^{1/4}+t^{1/2}).
$$
In 2002 Takahashi \cite{T. Miyakawa} studied the 
existence and space-time decay rates of global strong 
solutions of the Cauchy problem for the
Navier-Stokes equations in the weighted 
$L^\infty(\mathbb R^d, (1+|x|)^\beta{\rm dx})$ spaces. 
Takahashi showed that if $u_0$ satisfies
\begin{equation}\label{eq1}
|(e^{t\Delta}u_0)(x)| < \delta(1+|x|)^{-\beta}, 
|(e^{t\Delta}u_0)(x)| < \delta(1+t)^{-\frac{\beta}{2}},
\end{equation}
with sufficiently small $\delta$, then NSE has a global mild solution $u$ such that
$$
|u(x,t)|  \leq C(1+|x|)^{-\beta}, |u(x,t)|   \leq C(1+t)^{-\frac{\beta}{2}}, 
$$
where $\beta$ is restricted by the condition $1\leq \beta \leq d+1$.\\
Takahashi also showed that if 
$$
|u_0(x)| \leq c(1+|x|)^{-\beta} \ \ \ \ \ \ \ {\rm for\ some}\ 0<\beta \leq d,
$$ 
then 
$$
|(e^{t\Delta}u_0)(x)| \leq c(1+|x|)^{-\beta}, 
\ \ |(e^{t\Delta}u_0)(x)| \leq c(1+t)^{\frac{-\beta}{2}}.
$$
In this paper, we discuss the existence and space-time decay 
rates of global strong solutions of the Cauchy problem for the
Navier-Stokes equations in the weighted $L^\infty(\mathbb 
R^d, |x|^\gamma{\rm dx})\cap L^\infty(\mathbb R^d, |x|^\beta{\rm dx})$ spaces.
The spaces $L^\infty(\mathbb R^d, |x|^\gamma{\rm dx})
\cap L^\infty(\mathbb R^d, |x|^\beta{\rm dx})$ are more general than the spaces 
$L^\infty(\mathbb R^d, (1+|x|)^\beta{\rm dx})$.\linebreak In particular, $L^\infty(\mathbb R^d, |x|^\gamma{\rm dx})\cap L^\infty(\mathbb R^d, |x|^\beta{\rm dx})$
= $L^\infty(\mathbb R^d, (1+|x|)^\beta{\rm dx})$ when $\gamma =0$, 
and so this result improves the previous one.\\
The content of this paper is as follows: in Section 2, we state our main
theorems after introducing some notations. In Section 3, 
we first prove the  some estimates concerning the heat
semigroup with the Helmholtz-Leray \linebreak projection and some auxiliary lemmas. 
Finally, in Section 4, 
we will give the proof of the main theorems.
\vskip 0,5cm
\centerline{\bf \S2. Statement of the results} 
\vskip 0,5cm
Now, for $T > 0$, we say that $u$ is a mild solution of NSE on $[0, T]$ 
corresponding to a divergence-free initial datum $u_0$ 
when $u$ solves the integral equation
$$
u = e^{t\Delta}u_0 - \int_{0}^{t} e^{(t-\tau) \Delta} \mathbb{P} 
\nabla  .\big(u(\tau,.)\otimes u(\tau,.)\big) \dif\tau.
$$
Above we have used the following notation: 
for a tensor $F = (F_{ij})$ 
we define the vector $\nabla.F$ by 
$(\nabla.F)_i = \sum_{j = 1}^d\partial_jF_{ij}$ 
and for two vectors $u$ and $v$, we define their tensor 
product $(u \otimes v)_{ij} = u_iv_j$. 
The operator $\mathbb{P}$ is the Helmholtz-Leray 
projection onto the divergence-free fields 
$$
(\mathbb{P}f)_j =  f_j + \sum_{1 \leq k \leq d} R_jR_kf_k, 
$$
where $R_j$ is the Riesz transforms defined as 
$$
R_j = \frac{\partial_j}{\sqrt{- \Delta}}\ \ {\rm i.e.} \ \  
\widehat{R_jg}(\xi) = \frac{i\xi_j}{|\xi|}\hat{g}(\xi).
$$
The heat kernel $e^{t\Delta}$ is defined as 
$$
e^{t\Delta}u(x) = ((4\pi t)^{-d/2}e^{-|.|^2/4t}*u)(x).
$$
For a space of functions defined on $\mathbb R^d$, say $E(\mathbb R^d)$, we will abbreviate it as $E$ and we do not distinguish between the vector-valued and scalar-value spaces of functions. Throughout the paper, we sometimes use the notation $A \lesssim B$ as an equivalent to $A \leq CB$ with a uniform constant $C$. The notation $A \simeq B$
means that $A \lesssim B$ and $B \lesssim A$. Let $\beta \geq 0$, we define the space $L^\infty(|x|^\beta{\rm dx}):= L^\infty(\mathbb R^d,|x|^\beta{\rm dx})$ which is made up by the measurable  functions $u$ such that  
$$
\big\|u\big\|_{L^\infty(|x|^\beta{\rm dx})}: = \underset{x \in \mathbb R^d}{\rm esssup}|x|^\beta|u(x)| < +\infty.
$$
Now we can state our result
\begin{Th}\label{th4} 
Assume that $d\geq 1$, and $0\leq \gamma \leq 1 \leq \beta <d$. Then for all $f \in  L^\infty(|x|^\gamma{\rm dx}) \cap L^\infty(|x|^\beta{\rm dx})$ we have
\begin{gather*}
\underset{x \in \mathbb R^d, t>0}{\rm sup}\big(|x|^{\tilde \gamma}t^{\frac{1}{2}(\gamma-  \tilde \gamma)} + |x|^\alpha t^{\frac{1}{2}(1-  \alpha)}+|x|^{\tilde\beta} t^{\frac{1}{2}(\beta-  \tilde \beta)}\big)|e^{t\Delta}f| \\ \lesssim  \big\|f\big\|_{L^\infty(|x|^\gamma{\rm dx})}+\big\|f\big\|_{L^\infty(|x|^\beta{\rm dx})}
\end{gather*}
 for $0 \leq \tilde\gamma \leq \gamma, 0\leq \alpha \leq 1$, and $0 \leq \tilde \beta \leq \beta$.
\end{Th}

\begin{Th}\label{th2}   Let $0\leq \gamma \leq 1\leq \beta <d$ be fixed, then for all $\tilde \gamma, \alpha$, and $\tilde \beta$ satisfying 
$$
0\leq\tilde \gamma \leq\gamma, \tilde \beta \geq 0, \beta -2 < \tilde \beta \leq \beta, 0 < \alpha <1, \  and\ \beta -\tilde \beta -1< \alpha < d-\tilde\beta,
$$ 
there exists a positive constant $\delta_{\gamma,\tilde \gamma,\alpha,\beta,\tilde \beta,d}$ such that for all
 $u_0 \in L^\infty(|x|^\gamma{\rm d}x) \cap L^\infty(|x|^\beta{\rm d}x)$  with ${\rm div}(u_0) = 0$ satisfying
\begin{equation}\label{eq4}
\underset{x \in \mathbb R^d, t>0}{\rm sup}\big(|x|^{\tilde \gamma}t^{\frac{1}{2}(\gamma-  \tilde \gamma)} + |x|^\alpha t^{\frac{1}{2}(1-  \alpha)}+|x|^{\tilde\beta} t^{\frac{1}{2}(\beta-  \tilde \beta)}\big)|e^{t\Delta}u_0|
\leq \delta_{\gamma,\tilde \gamma,\alpha,\beta,\tilde \beta,d},
\end{equation}
NSE has a global mild solution $u$ on $(0,\infty) \times \mathbb R^d$ such that
\begin{equation}\label{eq5}
\underset{x \in \mathbb R^d, t>0}{\rm sup}\big(|x|^{\gamma}+t^{\frac{\gamma}{2}} +|x|^\beta +t^{\frac{\beta}{2}}\big)|u(x,t)| <+\infty.
\end{equation}
\end{Th}
\begin{Not} Our result improves the previous result for $L^\infty(\mathbb R^d, (1+|x|)^\beta{\rm dx})$. This space,
studied in  \cite{T. Miyakawa}, is a particular case of the 
space $L^\infty(|x|^\gamma{\rm d}x) \cap L^\infty(|x|^\beta{\rm d}x)$ when $\gamma=0$. 
Furthermore, we prove that Takahashi's result holds true under a much weaker condition on the initial data. Indeed, from Lemma \ref{lem4} and Theorem \ref{th4}, it is easily seen that the condition \eqref{eq4} of \linebreak Theorem \ref{th2} is weaker than the condition \eqref{eq1}.
\end{Not}
\begin{Not}We invoke Theorem \ref{th4} to deduce that  if $u_0 \in L^\infty(|x|^\gamma{\rm dx}) \cap L^\infty(|x|^\beta{\rm dx})$ and $\big\|u_0\big\|_{L^\infty(|x|^\gamma{\rm dx})}+ \big\|u_0\big\|_{L^\infty(|x|^\beta{\rm dx})}$ is small enough then the \linebreak condition \eqref{eq4} of Theorem \ref{th2} is valid.
\end{Not}
\begin{Th}\label{th2a}   Let $1\leq \beta <d$ be fixed, then for all $\alpha$ satisfying $0<\alpha< 1$, there exists a positive constant $\delta_{\alpha,d}$ such that for all
 $u_0 \in L^\infty(|x|{\rm d}x) \cap L^\infty(|x|^\beta{\rm d}x)$  with ${\rm div}(u_0) = 0$ satisfying
\begin{equation}\label{eq4a}
\underset{x \in \mathbb R^d, t>0}{\rm sup}|x|^\alpha t^{\frac{1}{2}(1-  \alpha)}|e^{t\Delta}u_0|
\leq \delta_{\alpha,d},
\end{equation}
NSE has a global mild solution $u$ on $(0,\infty) \times \mathbb R^d$ such that
$$
\underset{x \in \mathbb R^d, t>0}{\rm sup}\big(|x|+t^{\frac{1}{2}}\big)|u(x,t)| <+\infty 
$$
and
$$
\underset{x \in \mathbb R^d, 0<t<T}{\rm sup}|x|^\beta|u(x,t)| <+\infty,\  for\ all\ T \in (0, \infty).
$$
\end{Th}
\begin{Not}We invoke Theorem \ref{th4} to deduce that  if $u_0 \in L^\infty(|x|{\rm dx})$ and $\big\|u_0\big\|_{L^\infty(|x|{\rm dx})}$ is small enough then the condition \eqref{eq4a} of Theorem \ref{th2a} is valid.
\end{Not}
\vskip 0.5cm
\centerline{\bf \S3. Some auxiliary results} 
\vskip 0.5cm
In this section we establish some auxiliary lemmas. 
We first prove a version of Young's inequality type 
for convolutions in $L^\infty(|x|^\beta{\rm dx})$ spaces. 
\begin{Lem}\label{lem1} 
Assume that $d\geq 1, 0<\alpha<d, 0<\beta<d$ and $\alpha +\beta >d$. 
Then for all  $f \in L^\infty(|x|^\alpha{\rm dx})$ 
and for all  $g \in L^\infty(|x|^\beta{\rm dx})$ 
 we have
$$
\big\|f*g\big\|_{L^\infty(|x|^{\alpha+\beta-d}{\rm dx})} 
\lesssim \big\|f\big\|_{L^\infty(|x|^{\alpha}{\rm dx})}
\big\|g\big\|_{L^\infty(|x|^{\beta}{\rm dx})}.
$$
\end{Lem}
\textbf{Proof}.  Since $f*g$ is bilinear on 
$L^\infty(|x|^\alpha{\rm dx}) \times L^\infty(|x|^\beta{\rm dx})$, 
we may assume $\big\|f\big\|_{L^\infty(|x|^\alpha{\rm dx})} 
=\big\|g\big\|_{L^\infty(|x|^\beta{\rm dx})}=1$. We have
\begin{gather*}
(f*g)(x) =\int_{\mathbb R^d}f(x-y)g(y){\rm d}y 
= \int_{|y|<\frac{|x|}{2}} + \int_{\frac{|x|}{2}<|y|<\frac{3|x|}{2}} +\int_{|y|>\frac{3|x|}{2}}
=I_1+I_2+I_3.
\end{gather*}
From 
$$
|f(x)| \leq |x|^{-\alpha},\ {\rm and}\ |g(x)| \leq |x|^{-\beta},
$$
we get
\begin{gather*}
|I_1| \leq \int_{|y|<\frac{|x|}{2}}\frac{{\rm d}y}{|x-y|^\alpha |y|^\beta} \leq
\frac{2^\alpha}{|x|^\alpha}\int_{|y|<\frac{|x|}{2}}\frac{{\rm d}y}{ |y|^\beta}\simeq 
\frac{1}{|x|^{\alpha+\beta-d}}.
\end{gather*}
\begin{gather*}
|I_2| \leq \int_{\frac{|x|}{2}<|y|<\frac{3|x|}{2}}\frac{{\rm d}y}{|x-y|^\alpha |y|^\beta} \leq
\frac{2^\beta}{|x|^\beta}\int_{|y|<\frac{5|x|}{2}}\frac{{\rm d}y}{ |y|^\alpha}\simeq 
\frac{1}{|x|^{\alpha+\beta-d}}.
\end{gather*}
\begin{gather*}
|I_3| \leq \int_{|y|>\frac{3|x|}{2}}\frac{{\rm d}y}{|x-y|^\alpha |y|^\beta} \leq 
3^\alpha\int_{|y|>\frac{3|x|}{2}}\frac{{\rm d}y}{|y|^\alpha |y|^\beta}\simeq 
\frac{1}{|x|^{\alpha+\beta-d}}.
\end{gather*}
We thus obtain
$$
|(f*g)(x)| \lesssim \frac{1}{|x|^{\alpha+\beta-d}}.
$$
The proof Lemma \ref{lem1} is complete. \qed\\
We now deduce the $L^\infty(|x|^\gamma{\rm dx})-L^\infty(|x|^\beta{\rm dx})$ 
estimate for the heat semigroup.
\begin{Lem}\label{lem2} 
Assume that $d\geq 1$ and  $0\leq \gamma \leq \beta <d$. 
Then for all \linebreak $f \in L^\infty(|x|^\beta{\rm dx})$ we have
\begin{equation}\label{eq6}
\big\|e^{t\Delta}f\big\|_{L^\infty(|x|^\gamma{\rm dx})} 
\lesssim t^{-\frac{1}{2}(\beta-\gamma)}
\big\|f\big\|_{L^\infty(|x|^\beta{\rm dx})},\ for\ t>0.
\end{equation}
\end{Lem}
\textbf{Proof}. We have 
$$
(e^{t\Delta}f)(x)=\int_{\mathbb R^d}\frac{1}{t^{d/2}}
E\big(\frac{x-y}{\sqrt{t}}\big)f(y){\rm d}y,\ {\rm where}\ E(x) 
= (4\pi )^{-\frac{d}{2}}e^{-\frac{|x|^2}{4}}. 
$$
Recall the simate
\begin{equation}\label{eq7}
t^{-\frac{d}{2}}e^{-\frac{|x|^2}{4t}} 
\lesssim |x|^{-\alpha}t^{-\frac{1}{2}(d-\alpha)}, \ {\rm for}\  0 \leq \alpha \leq d.
\end{equation}
We first consider the case $0<\gamma <\beta$. 
From the inequality \eqref{eq7} and Lemma \ref{lem1}, we have
$$
|(e^{t\Delta}f)(x)| \lesssim \int_{\mathbb R^d}
\frac{\big\|f\big\|_{L^\infty(|x|^\beta{\rm dx})}}
{t^{\frac{1}{2}(\beta-\gamma)}|x-y|^{\gamma+d-\beta}
|y|^{\beta}}{\rm d}y \lesssim 
t^{-\frac{1}{2}(\beta-\gamma)}|x|^{-\gamma}
\big\|f\big\|_{L^\infty(|x|^\beta{\rm dx})}.
$$
This proves \eqref{eq6}.\\
We consider the case $0=\gamma<\beta$. 
Applying Proposition 2.4 $(b)$ in (\cite{P. G. Lemarie-Rieusset 2002}, pp. 20) 
and note that
$|x|^{-\beta} \in L^{\frac{d}{\beta},\infty}$
$$
|e^{t\Delta}f (x)| \lesssim t^{-\frac{d}{2}}\big\|E\big(\frac{.}
{\sqrt t}\big)\big\|_{L^{\frac{d}{d-\beta},1}}\big\|f\big\|_{L^{\frac{d}{\beta},\infty}} \lesssim
t^{-\frac{\beta}{2}}\big\|E\big\|_{L^{\frac{d}{d-\beta},1}}\big\|f\big\|_{L^\infty(|x|^\beta{\rm dx})}.
$$
This proves \eqref{eq6}.\\
Suppose finally that $0 \leq \gamma =\beta$. We have
$$
\int_{\mathbb R^d}\frac{1}{t^{d/2}}E\big(\frac{x-y}
{\sqrt{t}}\big)f(y){\rm d}y = \int_{|y|<\frac{|x|}{2}} 
+ \int_{|y|>\frac{|x|}{2}} = I_1 +I_2.
$$
From the inequality \eqref{eq7}, we have
\begin{gather*}
|I_1| \lesssim \big\|f\big\|_{L^\infty(|x|^\beta{\rm dx})} 
\int_{|y|<\frac{|x|}{2}}|x-y|^{-d}|y|^{-\beta}{\rm d}y  \leq \notag \\
 \big\|f\big\|_{L^\infty(|x|^\beta{\rm dx})}\big(\frac{|x|}{2}\big)^{-d}
\int_{|y|<\frac{|x|}{2}}|y|^{-\beta}{\rm d}y 
\simeq  \big\|f\big\|_{L^\infty(|x|^\beta{\rm dx})} |x|^{-\beta}.
\end{gather*}
\begin{gather*}
|I_2| \leq  \big\|f\big\|_{L^\infty(|x|^\beta{\rm dx})}\int_{|y|>\frac{|x|}{2}}
\frac{1}{t^{d/2}}E\big(\frac{x-y}{\sqrt{t}}\big)|y|^{-\beta}{\rm d}y  \leq \notag\\
 \big\|f\big\|_{L^\infty(|x|^\beta{\rm dx})}\big(\frac{|x|}{2}\big)^{-\beta}
\int_{y \in \mathbb R^d}\frac{1}{t^{d/2}}E\big(\frac{y}{\sqrt{t}}\big){\rm d}y= 
C \big\|f\big\|_{L^\infty(|x|^\beta{\rm dx})}|x|^{-\beta}, \notag
\end{gather*}
where
$$
C =2^\beta\int_{y \in \mathbb R^d}E(y){\rm d}y<+\infty.
$$
Therefore,
$$
|e^{t\Delta}f (x)| \lesssim  \big\|f\big\|_{L^\infty(|x|^\beta{\rm dx})}|x|^{-\beta}.
$$
The proof of Lemma \ref{lem2} is complete. \qed \\

We now deduce the $L^\infty(|x|^\gamma{\rm dx})-L^\infty(|x|^\beta{\rm dx})$ 
estimate for the operator $e^{t\Delta} \mathbb{P}\nabla$.
As shown in \cite{P. G. Lemarie-Rieusset 2002}, the kernel function $F_t$ of  $e^{t\Delta} \mathbb{P}\nabla$ satisfies the following \linebreak
 inequalities
\begin{equation}\label{eq8}
F_t(x) = t^{-\frac{d+1}{2}}F\big(\frac{x}{\sqrt{t}}\big), |F(x)| \lesssim \frac{1}{(1+|x|)^{d+1}},
\end{equation}
\begin{equation}\label{eq9}
|F_t(x)| \lesssim |x|^{-\alpha}t^{-\frac{1}{2}(d+1-\alpha)}, \ {\rm for}\  0 \leq \alpha \leq d+1.
\end{equation}
By using the inequalities \eqref{eq8} and \eqref{eq9} and 
arguing as in the proof of \linebreak Lemma \ref{lem2},  
we can easily prove the following lemma.
\begin{Lem}\label{lem3} 
Assume that $d\geq 1$ and  $0\leq \gamma \leq \beta <d$. Then for all \linebreak $f \in L^\infty(|x|^\beta{\rm dx})$ we have
$$
\big\|e^{t\Delta} \mathbb{P}\nabla.f\big\|_{L^\infty(|x|^\gamma{\rm dx})} \lesssim t^{-\frac{1}{2}(\beta+1-\gamma)}\big\|f\big\|_{L^\infty(|x|^\beta{\rm dx})},\ {\rm for}\ t>0.
$$
\end{Lem}
\begin{Lem}\label{lem4}
Let $0\leq \gamma < \beta \leq d$. Assume that $f \in \mathcal{S}'(\mathbb R^d)$ and satisfies the following inequality
\begin{equation}\label{eq10}
\underset{x \in \mathbb R^d, t>0}{\rm sup}\big(|x|^{\gamma}+|x|^{\beta}\big)|(e^{t\Delta}f)(x)|  =C <+\infty,
\end{equation}
then
$$
f \in L^\infty(|x|^\gamma{\rm dx})\cap L^\infty(|x|^\beta{\rm dx})
$$
and
\begin{equation}\label{eq11}
\underset{x \in \mathbb R^d}{\rm esssup}\big(|x|^{\gamma}+|x|^{\beta}\big)|f(x)|  \leq C.
\end{equation}
\end{Lem}
\textbf{Proof}. Since $\frac{1}{|x|^{\gamma}+|x|^{\beta}} \in L^{\frac{d}{\beta},\infty}\cap  L^{\frac{d}{\gamma},\infty}$ and $L^{\frac{d}{\beta},\infty}\cap  L^{\frac{d}{\gamma},\infty}\subset L^q$ for all $q$ satisfying $\frac{d}{\beta} < q < \frac{d}{\gamma}$, it follows that $e^{t\Delta}f \in L^{\infty}(0,\infty;L^q)$ for all $q \in \big(\frac{d}{\beta}, \frac{d}{\gamma}\big)$, by a compactness theorem in Banach space, there exists  a sequence $t_k$ which converges to $0$ 
such that $e^{t_k\Delta}f $  converges weakly to $f'$ in $L^q$ with $f' \in L^q$. Since $e^{t\Delta}$ 
is a continuous semigroup on $\mathcal{S}'(\mathbb R^d)$, it follows that $f=f' \in L^q$. Since $e^{t\Delta}$ 
is a continuous semigroup on $L^q(\mathbb R^d), (1\leq q <\infty)$, we get
$$
\underset{k \rightarrow \infty}{\rm lim}\big\|e^{t_k\Delta}f -f\big\|_{L^q}=0,\ {\rm for}\ q\in \big(\frac{d}{\beta}, \frac{d}{\gamma}\big).
$$
Therefore, there exists a subsequence $t_{k_j}$ of the sequence $t_k$ such that
\begin{equation}\label{eq12}
\underset{j \rightarrow \infty}{\rm lim}(e^{t_{k_j}\Delta}f)(x)=f(x)\ {\rm for\ almost\ everywhere}\ x \in \mathbb{R}^d.
\end{equation}
The inequality \eqref{eq11} is deduced from equalities \eqref{eq10} and \eqref{eq12}. \qed 
\begin{Not} We invoke Lemma \ref{lem4} for $\gamma =0$ and Lemma \ref{lem2} for $\gamma =\beta$  to deduce that the condition \eqref{eq1} of Takahashi on the initial data  
is equivalent to the condition 
$$
\big\|u_0\big\|_{L^\infty((1+|x|)^\beta{\rm dx})} \leq \delta.
$$
\end{Not}
\begin{Lem}\label{lem5} Let $\gamma, \theta \in \mathbb R$ and $t>0$, then\\
{\rm (a)} If $\theta< 1$ then
$$
\int^{\frac{t}{2}}_0(t-\tau)^{-\gamma} \tau^{-\theta}{\rm d}\tau = C t^{1-\gamma - \theta},\  where\  
C = \int^{\frac{1}{2}}_0(1-\tau)^{-\gamma} \tau^{-\theta}{\rm d}\tau < \infty.
$$
{\rm (b)} If $\gamma < 1$ then
$$
\int^{t}_{\frac{t}{2}}(t-\tau)^{-\gamma} \tau^{-\theta}{\rm d}\tau = C t^{1-\gamma - \theta},\  where\  
C = \int^1_{\frac{1}{2}}(1-\tau)^{-\gamma} \tau^{-\theta}{\rm d}\tau < \infty.
$$
{\rm (c)} If $\gamma < 1$ and  $\theta < 1$ then
$$
\int^{t}_0(t-\tau)^{-\gamma} \tau^{-\theta}{\rm d}\tau = C t^{1-\gamma - \theta},\  where\  
C = \int^1_0(1-\tau)^{-\gamma} \tau^{-\theta}{\rm d}\tau < \infty.
$$
\end{Lem}
The proof of this lemma is elementary and may be omitted.\qed \\
Let us recall the following result on solutions of a quadratic
equation in Banach spaces (Theorem 22.4 in 
\cite{P. G. Lemarie-Rieusset 2002}, p. 227).
\begin{Th}\label{th3}
Let $E$ be a Banach space, and $B: E \times E \rightarrow  E$ 
be a continuous bilinear map such that there exists $\eta > 0$ so that
$$
\|B(x, y)\| \leq \eta \|x\| \|y\|,
$$
for all x and y in $E$. Then for any fixed $y \in E$ 
such that $\|y\| \leq \frac{1}{4\eta}$, the equation $x = y - B(x,x)$ 
has a unique solution  $\overline{x} \in E$ satisfying 
$\|\overline{x}\| \leq \frac{1}{2\eta}$.
\end{Th}
\vskip 0.5cm
 \centerline{\bf\S4. Proofs of Theorems \ref{th4}, \ref{th2}, and \ref{th2a}}
\vskip 0.5cm
In this section we will give the proofs of  
Theorems  \ref{th4},  \ref{th2}, and \ref{th2a}. We now need eight more lemmas. 
In order to proceed, we define an auxiliary space $K^\beta_{\alpha,T}$. 
Let $\alpha, \beta$, and $T$ be such that $0 \leq \alpha \leq \beta<d, 0<T\leq +\infty$, 
we define the auxiliary space $K^\beta_{\alpha,T}$ 
which is made up by the measurable functions $u(t,x)$ such that 
\begin{gather*}
\underset{x \in \mathbb R^d, 0<t<T}{\rm esssup}
|x|^{\alpha}t^{\frac{1}{2}(\beta-\alpha)}|u(x,t)| < +\infty.
\end{gather*}
The auxiliary space $K^\beta_{\alpha,T}$ is equipped with the norm
$$
\big\|u\big\|_{K^\beta_{\alpha,T}}:= 
\underset{x \in \mathbb R^d, 0<t<T}{\rm esssup}
|x|^{\alpha}t^{\frac{1}{2}(\beta-\alpha)}|u(x,t)|.
$$
We rewrite Lemma  \ref{lem2}  as follows
\begin{Lem}\label{lem6} 
Assume that $d\geq 1$ and  $0\leq \alpha \leq \beta <d$. 
Then for all \linebreak  $f \in L^\infty(|x|^\beta{\rm dx})$ we have 
$e^{t\Delta}f \in K^\beta_{\alpha,T}$ and 
$\big\|e^{t\Delta}f\big\|_{K^\beta_{\alpha,T}} 
\leq C\big\|f\big\|_{L^\infty(|x|^\beta{\rm dx})}$, where C is a positive constant independent of T.
\end{Lem}
\begin{Lem}\label{lem7} 
Assume that $d\geq 1$ and  $0\leq \alpha \leq \beta <d$. 
Then 
$$K^\beta_{\alpha,T} \subset K^\beta_{\beta,T} \cap K^\beta_{0,T}.$$
\end{Lem}
The proof of this lemma is elementary and may be omitted.\qed \\
\begin{Lem}\label{lem7a} 
Assume that $d\geq 1,T<+\infty$, and  $0\leq \alpha \leq \beta \leq \tilde\beta <d$. 
Then $K^\beta_{\alpha,T}\subset  K^{\tilde\beta}_{\alpha,T}$.
\end{Lem}
The proof of this lemma is elementary and may be omitted.\qed \\
In the following lemmas a particular attention 
will be devoted to the study of the bilinear operator 
$B(u, v)(t)$ defined by 
\begin{equation}\label{eq13}
B(u, v)(t) = \int_{0}^{t} e^{(t-\tau ) \Delta} \mathbb{P} 
\nabla.\big(u(\tau)\otimes v(\tau)\big) \dif\tau.
\end{equation}
\begin{Lem}\label{lem8}
Let $\beta, \tilde \beta, \hat \beta$, and $\alpha$ be such that
\begin{gather*}
0 \leq \beta <d, \tilde \beta > \beta -2, 0\leq \tilde \beta \leq \beta, 0<\alpha<1, \beta -\tilde \beta-1 <\alpha < d-\tilde \beta ,\\
0 \leq \hat \beta \leq \beta,\  and\ \alpha +\tilde \beta-1 < \hat \beta \leq \alpha +\tilde \beta.
\end{gather*}
Then the bilinear operator $B$ is continuous from $K^1_{\alpha,T} \times K^\beta_{\tilde \beta,T}$ into $K^\beta_{\hat \beta,T}$
and the following inequality holds 
\begin{equation}\label{eq14}
\big\|B(u, v)\big\|_{K^\beta_{\hat \beta,T}} \leq C\big\|u\big\|_{K^1_{\alpha,T}}
\big\|v\big\|_{K^\beta_{\tilde \beta,T}},
\end{equation}
where C is a positive constant independent of T.
\end{Lem}
\textbf{Proof}. Since $B(.,.)$ is bilinear on 
$K^1_{\alpha,T} \times K^\beta_{\tilde \beta,T}$, 
we may assume  $\big\|u\big\|_{K^1_{\alpha,T}} =\big\|v\big\|_{K^\beta_{\tilde \beta,T}}=1$.
From 
$$
|(u\otimes v)| \leq |y|^{-(\alpha+\tilde \beta)}t^{-\frac{1}{2}(1-\alpha+\beta-\tilde \beta)}, 
$$
by using Lemma \ref{lem3}, we have
$$
\big|e^{(t-\tau ) \Delta} \mathbb{P}\nabla.\big(u\otimes v\big)\big| 
\lesssim |x|^{-\hat \beta}\frac{1}{(t-s)^{\frac{1}{2}
(1+\alpha+\tilde \beta-\hat \beta)}t^{\frac{1}{2}(1-\alpha+\beta-\tilde \beta)}}
$$
then applying Lemma \ref{lem5} (c), we get
$$
|B(u,v)| \lesssim |x|^{-\hat \beta}\int^t_0\frac{1}{(t-s)^{\frac{1}{2}
(1+\alpha+\tilde \beta-\hat \beta)}t^{\frac{1}{2}(1-\alpha+\beta-\tilde \beta)}}
{\rm d}s \simeq |x|^{-\hat \beta}t^{-\frac{1}{2}(\beta-\hat \beta)}.
$$
This proves Lemma \ref{lem8}.\\
Note that since $\alpha > \beta -\tilde \beta -1$ and $\hat \beta > \alpha+\tilde \beta-1$, 
it follows that the conditions $\frac{1-\alpha +\beta -\tilde \beta}{2} < 1$ and 
$\frac{1+\alpha +\tilde \beta -\hat \beta}{2} < 1$ are valid. So we can apply Lemma \ref{lem5} (c). \qed
\begin{Lem}\label{lem9} Assume that NSE has a mild solution $u \in K^1_{\tilde \alpha,T}$ for some \linebreak $\tilde \alpha \in (0,1)$ 
with initial data $u_0 \in L^\infty({|x|{\rm d}x})$ then $u \in K^1_{\alpha,T}$ for all $\alpha \in [0,1]$.
\end{Lem}
\textbf{Proof}. From $u = e^{t\Delta}u_0 + B(u,u)$, applying Lemmas \ref{lem6} and \ref{lem8} 
with $\beta =1$ and $\alpha =\tilde \beta =\tilde \alpha$, we get $u \in K^1_{\hat \beta,T}$
 for all $\hat \beta \in \big(\tilde \alpha-(1-\tilde \alpha), 2\tilde \alpha\big) \cap[0,1]$.
 Applying again Lemmas \ref{lem6} and \ref{lem8}  with $\beta =1, \alpha =\tilde \alpha$, and $\tilde \beta \in \big(\tilde \alpha-(1-\tilde \alpha), 2\tilde \alpha\big) \cap[0,1]$
 to get $u \in K^1_{\hat \beta,T}$ for all $\hat \beta \in \big(\tilde \alpha-2(1-\tilde \alpha), 3\tilde \alpha\big) \cap[0,1]$.
 By induction, we get $u\in K^1_{\hat \beta,T}$ for all $\hat \beta \in \big(\tilde \alpha-n(1-\tilde \alpha), (n+1)\tilde \alpha\big) \cap[0,1]$ with $n\in \mathbb N$. 
Since $\tilde\alpha \in (0,1)$, it follows that there exists sufficiently large $n$ satisfying
$$
\big(\tilde \alpha-n(1-\tilde \alpha), (n+1)\tilde \alpha\big) \supset [0,1].
$$
This proves Lemma \ref{lem9}. \qed
\begin{Lem}\label{lem10}Let $\beta$ be a fixed number in the interval $[0, d)$. Assume that NSE has a mild solution
 $u \in \underset{\alpha \in [0,1]}{\cap}K^1_{\alpha,T}\cap K^{\beta}_{\tilde \beta,T}$ for some $\tilde \beta \in [0,\beta] \cap (\beta -2, \beta]$ 
with initial data $u_0 \in L^\infty({|x|^\beta{\rm d}x})$, then $u \in K^\beta_{\hat \beta,T}$ for all $\hat\beta \in [0,\beta]\cap(\tilde \beta-1,\beta]$.
\end{Lem}
\textbf{Proof}. We first prove that $u\in K^\beta_{\hat \beta,T}$ for all $\hat \beta \in [0,\beta]\cap (\tilde \beta -1,\tilde \beta +1)$.\\
Let $\alpha_1$ and $\alpha_2$ be such that
$$
{\rm max}\{\beta -\tilde \beta-1,\hat \beta -\tilde \beta,0\} <\alpha_1<1
$$
and
$$
{\rm max}\{\hat \beta -\tilde \beta,0\} <\alpha_2<{\rm min}\{1,\hat\beta -\tilde\beta +1\}.
$$
We split the integral given in \eqref{eq13} into two parts coming from the subintervals $(0, \frac{t}{2})$ and $(\frac{t}{2}, t)$
\begin{gather*}
B(u, u)(t) = 
 \int_{0}^{\frac{t}{2}} e^{(t-\tau ) \Delta} \mathbb{P} 
\nabla.\big(u\otimes u\big) \dif\tau + \int_{\frac{t}{2}}^{t} e^{(t-\tau ) \Delta} \mathbb{P} 
\nabla.\big(u\otimes u\big) \dif\tau =I_1+I_2.
\end{gather*}
Since $u \in \underset{\alpha \in [0,1]}{\cap}K^1_{\alpha,T}$, it follows that
\begin{equation}\label{eq15}
|u(x,t)| \lesssim |x|^{-\alpha_1}t^{-\frac{1}{2}(1-\alpha_1)},
\end{equation}
\begin{equation}\label{eq16}
|u(x,t)| \lesssim |x|^{-\alpha_2}t^{-\frac{1}{2}(1-\alpha_2)},
\end{equation}
and since $u\in K^{\beta}_{\tilde \beta,T}$, it follows that
\begin{equation}\label{eq17}
|u(x,t)| \lesssim |x|^{-\tilde \beta}t^{-\frac{1}{2}(\beta-\tilde \beta)}.
\end{equation}
From the inequalities \eqref{eq15} and \eqref{eq17}, and Lemma \ref{lem3}, we get
$$
\big|e^{(t-\tau ) \Delta} \mathbb{P}\nabla.\big(u\otimes u\big)\big| \lesssim 
|x|^{-\hat \beta}\frac{1}{(t-s)^{\frac{1}{2}(1+\alpha_1+\tilde \beta-\hat \beta)}t^{\frac{1}{2}(1-\alpha_1+\beta-\tilde \beta)}}.
$$
Then applying Lemma \ref{lem5} (a), we have
\begin{equation}\label{eq18}
|I_1| \lesssim |x|^{-\hat \beta}\int^{\frac{t}{2}}_0\frac{1}{(t-s)^{\frac{1}{2}(1+\alpha_1+\tilde \beta-\hat \beta)}t^{\frac{1}{2}(1-\alpha_1+\beta-\tilde \beta)}}{\rm d}s \simeq |x|^{-\hat \beta}t^{-\frac{1}{2}(\beta-\hat \beta)}.
\end{equation}
From the inequalities \eqref{eq16} and \eqref{eq17}, and Lemma \ref{lem3}, we get
$$
\big|e^{(t-\tau ) \Delta} \mathbb{P}\nabla.\big(u\otimes u\big)\big| \lesssim 
|x|^{-\hat \beta}\frac{1}{(t-s)^{\frac{1}{2}(1+\alpha_2
+\tilde \beta-\hat \beta)}t^{\frac{1}{2}(1-\alpha_2+\beta-\tilde \beta)}}.
$$
Then applying Lemma \ref{lem5} (b), we have
\begin{equation}\label{eq19}
|I_2| \lesssim |x|^{-\hat \beta}\int^{t}_{\frac{t}{2}}\frac{1}{(t-s)^{\frac{1}{2}
(1+\alpha_2+\tilde \beta-\hat \beta)}t^{\frac{1}{2}
(1-\alpha_2+\beta-\tilde \beta)}}{\rm d}s 
\simeq |x|^{-\hat \beta}t^{-\frac{1}{2}(\beta-\hat \beta)}.
\end{equation}
From the inequalities \eqref{eq18} and \eqref{eq19}, 
we get $B(u,u)\in K^\beta_{\hat \beta,T}$, 
and from \linebreak $u =e^{t\Delta}u_0+B(u,u)$ and Lemma \ref{lem6}, 
we have $u\in K^\beta_{\hat \beta,T}$.  This proves the result.\\
We now prove  $u \in K^\beta_{\hat \beta,T}$ for all 
$\hat\beta \in [0,\beta] \cap(\tilde\beta-1,\beta]$. Indeed,  
if $\tilde \beta > \beta -1$ then $u\in K^\beta_{\hat \beta,T}$ for all 
$\hat \beta \in [0,\beta]\cap (\tilde \beta -1,\tilde\beta+1) = [0,\beta]\cap (\tilde \beta -1,\beta]$ and so the 
lemma is proved. In the case $\tilde \beta \leq \beta -1$,  in exactly the same way, since $u\in K^\beta_{\hat \beta,T}$ 
for all $\hat \beta \in [0,\beta]\cap (\tilde \beta -1,\tilde\beta+1)$,
 it follows that $u\in K^\beta_{\hat \beta,T}$ for all 
$\hat \beta \in [0,\beta]\cap (\tilde \beta -1,\tilde \beta+2)
= [0,\beta]\cap (\tilde \beta -1,\beta]$. Therefore the proof 
of Lemma \ref{lem10} is complete. \qed
\begin{Lem}\label{lem11}Assume that NSE has a mild solution
 $u \in \underset{\alpha \in [0,1]}{\cap}K^1_{\alpha,T}\cap \underset{\hat \beta\in [\tilde\beta,\beta]}{\cap}K^{\beta}_{\hat \beta,T}$ 
for some $\tilde \beta \in [0,\beta]$ with initial data 
$u_0 \in L^\infty({|x|^\beta{\rm d}x})$. Then 
$u \in K^\beta_{\hat \beta,T}$ for all $\hat\beta \in [0,\beta]$.
\end{Lem}
\textbf{Proof}. We first prove that $u\in K^\beta_{\hat \beta,T}$ 
for all $\hat \beta \in [0,\beta]\cap(\tilde \beta-1,\tilde \beta]$.\\
We split the integral given in \eqref{eq13} into two parts coming from the subintervals $(0, \frac{t}{2})$ and $(\frac{t}{2}, t)$
\begin{gather*}
B(u, u)(t) = 
 \int_{0}^{\frac{t}{2}} e^{(t-\tau ) \Delta} \mathbb{P} 
\nabla.\big(u\otimes u\big) \dif\tau + \int_{\frac{t}{2}}^{t} e^{(t-\tau ) \Delta} \mathbb{P} 
\nabla.\big(u\otimes u\big) \dif\tau =I_1+I_2.
\end{gather*}
Let $\alpha_1$ be such that $0<\alpha_1<1$. Since 
$u \in K^1_{\alpha_1,T}\cap K^{\beta}_{\beta,T}$, it follows that
\begin{equation}\label{eq20}
|u(x,t)| \lesssim |x|^{-\alpha_1}t^{-\frac{1}{2}(1-\alpha_1)},
\end{equation}
\begin{equation}\label{eq21}
|u(x,t)| \lesssim |x|^{-\beta}.
\end{equation}
From the inequalities \eqref{eq20} and \eqref{eq21}, 
and Lemma \ref{lem3}, we get
$$
\big|e^{(t-\tau ) \Delta} \mathbb{P}\nabla.\big(u\otimes u\big)\big| \lesssim 
|x|^{-\hat \beta}\frac{1}{(t-s)^{\frac{1}{2}
(1+\alpha_1+\beta-\hat \beta)}t^{\frac{1}{2}(1-\alpha_1)}}.
$$
Then applying Lemma \ref{lem5} (a), we have
\begin{equation}\label{eq22}
|I_1| \lesssim |x|^{-\hat \beta}\int^{\frac{t}{2}}_0
\frac{1}{(t-s)^{\frac{1}{2}(1+\alpha_1+\beta-\hat \beta)}
t^{\frac{1}{2}(1-\alpha_1)}}{\rm d}s \simeq |x|^{-\hat \beta}
t^{-\frac{1}{2}(\beta-\hat \beta)}.
\end{equation}
Since $u \in K^1_{0,T}\cap K^{\beta}_{\tilde\beta,T}$, it follows that
\begin{equation}\label{eq23}
|u(x,t)| \lesssim t^{-\frac{1}{2}}\ {\rm and}\  |u(x,t)| 
\lesssim |x|^{-\tilde \beta}t^{-\frac{1}{2}(\beta-\tilde\beta)}.
\end{equation}
From the inequality \eqref{eq23}, and Lemma \ref{lem3}, we get
$$
\big|e^{(t-\tau ) \Delta} \mathbb{P}\nabla.\big(u\otimes u\big)\big| \lesssim 
|x|^{-\hat \beta}\frac{1}{(t-s)^{\frac{1}{2}(1+\tilde\beta-\hat \beta)}
t^{\frac{1}{2}(1+\beta-\tilde\beta)}}.
$$
Then applying Lemma \ref{lem5} (b), we obtain
\begin{equation}\label{eq24}
|I_2| \lesssim |x|^{-\hat \beta}\int^{t}_{\frac{t}{2}}\frac{1}{(t-s)^{\frac{1}{2}
(1+\tilde\beta-\hat \beta)}t^{\frac{1}{2}(1+\beta-\tilde\beta)}}{\rm d}s 
\simeq |x|^{-\hat \beta}t^{-\frac{1}{2}(\beta-\hat \beta)}.
\end{equation}
From the inequalities \eqref{eq22} and \eqref{eq24}, 
we get $B(u,u)\in K^\beta_{\hat \beta,T}$. From \linebreak $u =e^{t\Delta}u_0+B(u,u)$ 
and Lemma \ref{lem6}, we deduce $u\in K^\beta_{\hat \beta,T}$. This proves the result. Therefore, we get 
$u\in K^\beta_{\hat \beta,T}$ for all $\hat \beta \in [0,\beta]\cap(\tilde \beta-1,\beta]$.\\
We now prove that $u\in K^\beta_{\hat \beta,T}$ for all $\hat \beta \in [0,\beta]$. 
Indeed, in exactly the same way, since $u\in K^\beta_{\hat \beta,T}$ 
for all $\hat \beta \in [0,\beta]\cap(\tilde \beta-1, \beta]$,
 it follows that $u\in K^\beta_{\hat \beta,T}$ for all 
$\hat \beta \in [0,\beta]\cap(\tilde \beta-2,\beta]$. By induction, 
we get $u\in K^\beta_{\hat \beta,T}$ for all 
$\hat \beta \in [0,\beta]\cap (\tilde \beta -n, \beta]$ with $n \in \mathbb N$. 
However, there exists a sufficiently large number 
$n$ satisfying $\tilde\beta-n<0$ and therefore $u\in K^\beta_{\hat \beta,T}$ for all 
$\hat \beta \in [0,\beta]$. \linebreak The proof of Lemma \ref{lem11} is complete. \qed
\begin{Lem}\label{lem12}  
Let $0\leq \beta <d$ be fixed, then for all $\alpha$ and $\tilde \beta$ satisfying
$$
\tilde \beta \geq 0, 0 < \alpha <1, \beta -2 < \tilde \beta \leq \beta,
\  and\ \beta -\tilde \beta -1< \alpha < d-\tilde\beta,
$$ 
there exists a positive constant $\delta_{\alpha,\beta,\tilde \beta,d}$ such that for all
 $u_0 \in L^\infty(|x|{\rm d}x) \cap L^\infty(|x|^\beta{\rm d}x)$  
with ${\rm div}(u_0) = 0$ satisfying
\begin{equation}\label{eq2}
\underset{x \in \mathbb R^d, t>0}{\rm sup}\big(|x|^{\alpha}
t^{\frac{1}{2}(1-  \alpha)} +|x|^{\tilde \beta} 
t^{\frac{1}{2}(\beta-  \tilde \beta)}\big)|e^{t\Delta}u_0|
\leq \delta_{\alpha,\beta,\tilde \beta,d},
\end{equation}
NSE has a global mild solution $u$ on $(0,\infty) \times \mathbb R^d$ such that
\begin{equation}\label{eq3}
\underset{x \in \mathbb R^d, t>0}{\rm sup}\big(|x|+t^{\frac{1}{2}} 
+|x|^\beta +t^{\frac{\beta}{2}}\big)|u(x,t)| <+\infty.
\end{equation}
\end{Lem}
\textbf{Proof}. Applying Lemma \ref{lem8}  we deduce that the bilinear operator $B$ is bounded from  $K^1_{\alpha,\infty} \times K^1_{\alpha,\infty}$ 
into $K^1_{\alpha, \infty}$ and from  $K^1_{\alpha,\infty}  \times K^\beta_{\tilde \beta,\infty}$ into $K^\beta_{\tilde \beta,\infty}$. Therefore, 
the  bilinear operator $B$ is bounded from 
$$
(K^1_{\alpha,\infty} \cap K^\beta_{\tilde \beta,\infty})\times (K^1_{\alpha,\infty} \cap K^\beta_{\tilde \beta,\infty}) \ {\rm into}\
(K^1_{\alpha,\infty} \cap K^\beta_{\tilde \beta,\infty}).
$$ 
where the space $K^1_{\alpha,\infty} \cap K^\beta_{\tilde \beta,\infty}$ is equipped with the norm
$$
\big\|u\big\|_{K^1_{\alpha,\infty} \cap K^\beta_{\tilde \beta,\infty}}: = {\rm max}\{\big\|u\big\|_{K^1_{\alpha,\infty}},\big\|u\big\|_{K^\beta_{\tilde \beta,\infty}}\}.
$$
Applying Theorem \ref{th3} to the  bilinear operator $B$, we deduce that there exists a positive constant $\delta_{\alpha,\beta,\tilde \beta,d}$ 
such that for all $u_0 \in L^\infty({|x|{\rm d}x})\cap L^\infty({|x|^\beta{\rm d}x})$  with ${\rm div}(u_0) = 0$ satisfying 
$$
\big\|e^{t\Delta}u_0\big\|_{K^1_{\alpha,\infty} \cap K^\beta_{\tilde \beta,\infty}} \leq \delta_{\alpha,\beta,\tilde \beta,d},
$$
then NSE has a unique mild solution $u$ satisfying
$$
u \in K^1_{\alpha,\infty} \cap K^\beta_{\tilde \beta,\infty}.
$$
Applying Lemmas \ref{lem9}, \ref{lem10}, and \ref{lem11}, we get $u \in K^\beta_{\hat \beta,\infty}$ for all $\hat\beta \in [0,\beta]$ and $u \in K^1_{\alpha,\infty}$ for all $\alpha \in [0,1]$. The proof of Lemma \ref{lem12} is now complete. \qed 
\vskip 0.5cm
{\bf Proof of Theorem \ref{th4}}
\vskip 0.5cm
Since $|x| \leq |x|^\gamma +|x|^\beta$, it follows that 
$$
\big\|f\big\|_{L^\infty}(|x|{\rm d}x) \leq \big\|f\big\|_{L^\infty}(|x|^\gamma{\rm d}x)+\big\|f\big\|_{L^\infty}(|x|^\beta{\rm d}x).
$$
From Lemma \ref{lem2} we have
$$
|x|^\alpha t^{\frac{1}{2}(1-  \alpha)}|e^{t\Delta}u_0| \lesssim \big\|f\big\|_{L^\infty}(|x|{\rm d}x) \leq \big\|f\big\|_{L^\infty}(|x|^\gamma{\rm d}x)+\big\|f\big\|_{L^\infty}(|x|^\beta{\rm d}x),
$$
$$
|x|^{\tilde \gamma}t^{\frac{1}{2}(\gamma-  \tilde \gamma)}|e^{t\Delta}u_0|  \lesssim \big\|f\big\|_{L^\infty}(|x|^\gamma{\rm d}x),\ {\rm and}\ |x|^{\tilde \beta}t^{\frac{1}{2}(\beta-  \tilde \beta)}|e^{t\Delta}u_0| \lesssim  \big\|f\big\|_{L^\infty}(|x|^\beta{\rm d}x).
$$
This proves Theorem \ref{th4}.\qed 

\vskip 0.5cm
{\bf Proof of Theorem \ref{th2}}
\vskip 0.5cm
Since $L^\infty({|x|{\rm d}x})\subset L^\infty({|x|^\gamma{\rm d}x})\cap L^\infty({|x|^\beta{\rm d}x})$, it follows that  $u_0 \in L^\infty({|x|{\rm d}x}) $. Applying Lemma \ref{lem12} then there exists 
a positive constant $\delta_{\alpha,\beta, \tilde \beta,d}$ such that if
$$
\underset{x \in \mathbb R^d, t>0}{\rm sup}\big(|x|^{\alpha}t^{\frac{1}{2}(1-  \alpha)} +|x|^{\tilde \beta} t^{\frac{1}{2}(\beta-  \tilde \beta)}\big)|e^{t\Delta}u_0| \leq \delta_{\alpha,\beta, \tilde \beta,d},
$$
NSE has a global mild solution $u$ on $(0,\infty) \times \mathbb R^d$ such that
$$
\underset{x \in \mathbb R^d, t>0}{\rm sup}\big(|x|+t^{\frac{1}{2}}+|x|^{\beta}+t^{\frac{\beta}{2}}\big)|u(x,t)| <+\infty.
$$
Applying Lemma \ref{lem12} for $\beta =\gamma$ then there exists 
a positive constant $\delta_{\alpha,\gamma, \tilde \gamma, d}$ such that if
$$
\underset{x \in \mathbb R^d, t>0}{\rm sup}\big(|x|^{\alpha}t^{\frac{1}{2}(1-  \alpha)}+|x|^{\tilde \gamma} t^{\frac{1}{2}(\gamma-  \tilde \gamma)}\big)|e^{t\Delta}u_0|\leq \delta_{\alpha,\gamma, \tilde \gamma, d},
$$
NSE has a global mild solution $u$ on $(0,\infty) \times \mathbb R^d$ such that
$$
\underset{x \in \mathbb R^d, t>0}{\rm sup}\big(|x| +t^{\frac{1}{2}}+|x|^{\gamma}+t^{\frac{\gamma}{2}}\big)|u(x,t)| <+\infty.
$$
Therefore, if $u_0$ satisfies the following inequality
$$
\underset{x \in \mathbb R^d, t>0}{\rm sup}\big(|x|^{\tilde \gamma}t^{\frac{1}{2}(\gamma-  \tilde \gamma)} + |x|^\alpha t^{\frac{1}{2}(1-  \alpha)}+|x|^{\tilde\beta} t^{\frac{1}{2}(\beta-  \tilde \beta)}\big)|e^{t\Delta}u_0| \leq
{\min}\{\delta_{\alpha,\beta, \tilde \beta,d},\delta_{\alpha,\gamma, \tilde \gamma, d}\}
$$
NSE has a global mild solution $u$ on $(0,\infty) \times \mathbb R^d$ such that  \eqref{eq5}. \\
The proof of Theorem \ref{th2} is complete.\qed
\vskip 0.5cm
{\bf Proof of Theorem \ref{th2a}}
\vskip 0.2cm 
Applying Lemma \ref{lem8}  we deduce that the bilinear operator $B$ is bounded from  $K^1_{\alpha,\infty} \times K^1_{\alpha,\infty}$ 
into $K^1_{\alpha, \infty}$. Applying Theorem \ref{th3} to the  bilinear operator $B$, we deduce that there exists a positive constant $\delta_{\alpha,d}$ 
such that for all $u_0 \in L^\infty({|x|{\rm d}x})$  with ${\rm div}(u_0) = 0$ satisfying 
$$
\big\|e^{t\Delta}u_0\big\|_{K^1_{\alpha,\infty}} \leq \delta_{\alpha,d},
$$
then NSE has a unique mild solution $u$ satisfying $u \in K^1_{\alpha,\infty}$. Applying \linebreak Lemma \ref{lem9} we have $u \in \underset{\alpha \in [0,1]}{\cap}K^1_{\alpha,\infty}$.\\
We prove that  $u \in K^{\beta}_{\beta, T}$ for all $T \in (0, \infty)$.
Indeed, let $\gamma$ be such that $\gamma \in [1,\beta]\cap(\alpha,\alpha+2)$. Applying Lemma \ref{lem7a} we have $u \in K^{\gamma}_{\alpha,T}$, then using Lemmas \ref{lem10} we get $u \in K^{\gamma}_{\gamma, T}$,  in exactly the same way, using again Lemmas \ref{lem10},  since $u\in K^{\gamma}_{\gamma,T}$ for $\gamma \in [1,\beta]\cap(\alpha,\alpha+2)$, it follows that $u\in K^{\gamma}_{\gamma,T}$ for $\gamma \in [1,\beta]\cap(\alpha,\alpha+4)$. By induction, we get $u\in K^{\gamma}_{\gamma,T}$ for $\gamma \in [1,\beta]\cap(\alpha,\alpha+2n)$. However, there exists a sufficiently large number 
$n$ satisfying $\alpha+2n >\beta$ and therefore  $u \in K^{\beta}_{\beta, T}$. The proof of Theorem  \ref{th2a} is complete. \qed
\vskip 0.5cm
{\bf Acknowledgments}. This research is funded by 
Vietnam National \linebreak Foundation for Science and Technology 
Development (NAFOSTED) under grant number  101.02-2014.50.

\end{document}